\documentclass[a4paper,12pt]{article}
\usepackage[T1]{fontenc}
\usepackage[latin1]{inputenc}
\usepackage[english]{babel}
\usepackage[dvipdfm,pdfstartview=FitH,colorlinks=true,linkcolor=Red4,urlcolor=Red4,%
filecolor=green,citecolor=red]{hyperref}
\usepackage{url}
\usepackage[dvips]{graphics}
\usepackage{amsmath}
\usepackage{amssymb}
\usepackage{amsopn}
\usepackage{amsbsy}
\usepackage{amstext}
\usepackage{latexsym}
\usepackage{amsfonts}
\usepackage{amsthm}
\usepackage{array}
\usepackage{epsfig}
\usepackage[usenames,x11names]{xcolor}
\usepackage[top=2cm,bottom=2cm,margin=2cm]{geometry}
\usepackage{mathrsfs} 
\title{On a conjecture of Erd\H{o}s concerning primitive sequences}
\author{\sc Bakir FARHI \\[1mm] 
\href{mailto:bakir.farhi@gmail.com}{bakir.farhi@gmail.com} \\[1mm]
\url{http://www.bakir-farhi.site}
} 
\date{}

\newtheorem{thm}{Theorem}
\newtheorem{prop}[thm]{Proposition}

\newtheorem{conj}[thm]{Conjecture}

\let\epsilon=\varepsilon
\def\A{{\mathscr{A}}}
\def\P{{\mathscr{P}}}

\def\EMdash{\leavevmode\hbox to 10.6mm{\vrule height .63ex depth -.59ex
    width 10mm\hfill}}

\def\erdos{Erd\H{o}s}

\begin{document}
\maketitle

\begin{abstract}
In this note, we propose a conjecture stating that some series involving primitive sequences are convergent. Then, we show (by a counterexample) that the analogue of a conjecture of \erdos{}, for those series, is false.
\end{abstract}


\bigskip

Throughout this note, we denote by $\P$ the sequence of the prime numbers. Further, for a given positive integer $n$, we respectively denote by $\omega(n)$ and $\Omega(n)$ the number of prime factors of $n$ and the number of prime factors of $n$ counted with multiplicity. For a given sequence of positive integers $\A$, the quantity defined by $d°(\A) := \max\{\Omega(a) , a \in \A\}$ is called \emph{the degree} of $\A$. Particularly, if $\Omega(a)$ is the same for any $a \in \A$, then $\A$ is called \emph{an homogeneous sequence}.

A sequence $\A$ of positive integers is called \emph{primitive} if no term of the sequence divides any other. In \cite{erd1}, Erd\H{o}s proved that for any primitive sequence $\A$ (with $A \neq \{1\}$), the series $\sum_{a \in \A}\frac{1}{a \log a}$ converges and its sum is bounded above by an absolute constant $C$; and in \cite{ez}, \erdos{} and Zhang showed that $C \leq 1.84$. This bound is later improved by Clark \cite{c} to $e^{\gamma}$ $(\simeq 1.78)$, where $\gamma$ is the Euler constant. Furthermore, in \cite{erd2}, \erdos{} asked if it is true that the sum $\sum_{a \in \A} \frac{1}{a \log a}$ (where $\A \neq \{1\}$ is a primitive sequence) reaches its maximum value at $\A = \P$. Some years later, \erdos{} and Zhang \cite{ez} conjectured an affirmative answer to the last question by proposing the following:

\medskip

\noindent\textbf{Conjecture (\erdos{}):}~\\
\emph{%
For any primitive sequence $\A \neq \{1\}$, we have:
$$
\sum_{a \in \A} \frac{1}{a \log a} \leq \sum_{p \in \P} \frac{1}{p \log p} .
$$
}

\medskip

To compare with Clark's upper bound, let us precise that $\sum_{p \in \P} \frac{1}{p \log p} \simeq 1.63$. In \cite{z1}, Zhang proved the \erdos{} conjecture for a primitive sequence $\A$ ($\A \neq \{1\}$) satisfying $d°(\A) \leq 4$ and in \cite{z2}, he proved it for the particular case of homogeneous sequences and for some other primitive sequences slightly more complicated. To our knowledge, these are the only significant results that were obtained in the direction of proving \erdos{}'s conjecture.

In this note, we propose a conjecture stating that some series involving primitive sequences are convergent. Although such series are close to the series $\sum_{a \in \A} \frac{1}{a \log a}$, we will show (by a counterexample) that the analogue of the above \erdos{} conjecture for those series is false. We conjecture the following:
\begin{conj}\label{conj1}
For any primitive sequence $\A \neq \{1\}$, the series
$$
\sum_{a \in \A} \frac{\omega(a)}{a \log{a}} ~~~~\text{and}~~~~ \sum_{a \in \A} \frac{\Omega(a)}{a \log a}
$$
are both convergent.
\end{conj} 

Now, we will prove that the analogue of the above conjecture of \erdos{} for the series considered in Conjecture \ref{conj1} is false. In other words, we prove that the sum $\sum_{a \in \A} \frac{\omega(a)}{a \log a}$ (where $\A$ runs on the set of all primitive sequences different from $\{1\}$) does not reach its maximum value at $\A = \P$. The same will obviously be true for the sum $\sum_{a \in \A} \frac{\Omega(a)}{a \log a}$. We have the following:
\begin{prop}\label{p1}
There exists a primitive sequence $\A \neq \{1\}$ such that:
$$
\sum_{a \in \A} \frac{\omega(a)}{a \log a} ~>~ \sum_{p \in \P} \frac{\omega(p)}{p\log p} = \sum_{p \in \P} \frac{1}{p \log p} .
$$
\end{prop} 

\begin{proof}
Set
\begin{align*}
\A_1 &:= \{p q ~|~ p , q \in \P ~;~ p \leq 14 \times 10^5 ~\text{and}~ q \leq 14 \times 10^5\} \\
\A_2 &:= \{r \in \P ~|~ r > 14 \times 10^5\} \\
\A &:= \A_1 \cup \A_2 .
\end{align*}
It is clear that $\A$ is a primitive sequence. The computer calculations give:
$$
\sum_{a \in \A_1}\frac{\omega(a)}{a \log a} ~=~ \sum_{%
\begin{subarray}{c}
p , q \in \P \\
p \leq 14 \times 10^5 \\
q \leq 14 \times 10^5
\end{subarray}} \frac{1}{p q (\log p + \log q)} ~=~ 1.5748\dots ~~~~\text{and}~~~~
\sum_{%
\begin{subarray}{c}
p \in \P \\
p \leq 14 \times 10^5
\end{subarray}
}\frac{1}{p \log p} ~=~ 1.5659\dots
$$
So, we have: 
$$
\sum_{a \in \A_1}\frac{\omega(a)}{a \log a} > \sum_{%
\begin{subarray}{c}
p \in \P \\
p \leq 14 \times 10^5
\end{subarray}}\frac{1}{p \log p} .
$$
By adding $\displaystyle \sum_{a \in \A_2}\frac{\omega(a)}{a \log a} = \sum_{%
\begin{subarray}{c}
p \in \P \\
p > 14 \times 10^5
\end{subarray}} \frac{1}{p \log p}$ to the two hand-sides of this inequality, we get:
$$
\sum_{a \in \A}\frac{\omega(a)}{a \log a} > \sum_{p \in \P}\frac{1}{p \log p} ,
$$
as required.
\end{proof}

\end{document}